\documentclass[12pt]{amsart}

\usepackage{amssymb, latexsym, amsthm}

\newtheorem{thm}{Theorem}[section] %the resolution couldbe[subsection]
\newtheorem{prop}[thm]{Proposition}
\newtheorem{Prop}[thm]{Proposition}
\newtheorem{lemma}[thm]{Lemma}
\newtheorem{cor}[thm]{Corollary}

\theoremstyle{definition} %
\newtheorem{rem}[thm]{Remark}
\newtheorem{defn}[thm]{Definition}

\newcommand\CQ{{\mathcal{Q}}}
\newcommand\Deltah{{\Delta_{H}}}
\newcommand\tDeltah{{\tilde\Delta_{H}}}
\newcommand\IQ{{I\!\CQ}}
\newcommand\fp{{\mathfrak p}}
\newcommand\PCH{{\rm PCH}}
\newcommand\card[1]{{\left|{#1}\right|}}

\DeclareMathOperator{\Hur}{{\CQ_{Hur}}}

\DeclareMathOperator{\Hurone}{{\CQ_{Hur}^{1}}}

\renewcommand\th[1]{{${#1}^{\rm{th}}$}}
\def\Q{{\mathbb {Q}}}
\def\R{{\mathbb {R}}}
\def\C{{\mathbb {C}}}
\def\F{{\mathbb {F}}}
\def\Z{{\mathbb {Z}}}
\def\O{{\mathcal{O}}}
\def\HC{{\mathcal H}}
\def\sub{\subseteq}
\def\s{\sigma}

\newcommand\subjectto{{\,|\ }}
\newcommand\M[1][d]{{\operatorname{M}_{#1}}}
\newcommand\opname[1]{{\operatorname{#1}}}
\newcommand\Ker{{\opname{Ker}}}

\newcommand\PSL[1][d]{{\operatorname{PSL}_{#1}}}
\newcommand\SL[1][d]{{\operatorname{SL}_{#1}}}
\newcommand\dimcol[2]{{[{#1}\!:\!{#2}]}}
\newcommand\Tref[1]{{Theorem~\ref{#1}}}

\newcommand\Cref[1]{{Corollary~\ref{#1}}}
\newcommand\eq[1]{{(\ref{#1})}}

\newcommand\abs[2][]{\left|{#2}\right|_{#1}} % (choose which absolute valueisreferred to)
\def\normali{{\lhd}} % triangle without the - sign: for ideals.
\newcommand\ideal[1]{{\left<{#1}\right>}}
\def\divides{{\,|\,}} % spacing: 3/18.
\newcommand{\tr}[1][]{\operatorname{tr}_{#1}} %
\newcommand{\Norm}[1][]{{\operatorname{N}_{#1}}}
\newcommand\tensor[1][{}]{{\otimes_{#1}}}
\def\isom{{\;\cong\;}} % Meant to generate the = sign, with ~ above it
\newcommand{\set}[1]{{\{#1\}}}
\def\hra{{\,\hookrightarrow\,}}
\def\lam{{\lambda}}
\newcommand\mul[1]{{#1^{\times}}} % The multiplicative group
\def\({\left(}
\def\){\right)}
\newcommand\Pref[1]{{Proposition \ref{#1}}}
\newcommand{\Tr}[1][]{\operatorname{Tr}_{#1}^{\phantom{I}}} %
\long\def\forget#1\forgotten{{(\!(\!(\!(\!(\ \tiny{#1}
})\!)\!)\!)\!)\\ }

\def\ie {{\it i.e.\ }}
\def\eg {{\it e.g.\ }}
\def\cf {\hbox{\it cf.\ }}
\newcommand\md[2][d]{{\mul{#2}/{\mul{#2}}^{#1}}}
\def\co{{\,{:}\,}}
\def\ra{{\rightarrow}}
\newcommand\ver[1]{\marginpar{\tiny Changed in Ver \VER}}

\newcommand\algint[2][]{\if!#1\relax O_{#2}
\else{O_{\!#2}^{\phantom{I}}}\fi}  \DeclareMathOperator{\syspi}{{\rm sys}\pi}

\DeclareMathOperator{\vol}{{\rm vol}}

\newcommand\SLform[1]{{#1^{1}}} % Was: #1_1^{\times}.
\def\Lat{{\SLform{\CQ}}} \def\genus{{\operatorname{\it g}}}
\def\quat{{D}} % D is the common notation for a division algebra.
\newcommand\dom[2]{{{#1}\backslash{#2}}} \def\minusset{{-}}
\newcommand\suchthat{{\,:\ \,}} \newcommand{\smat}[4]
{{\(\!\!\begin{array}{cc}{#1}\!&\!{#2}\\[-0.1cm]{#3}\!&\!{#4}\end{array}\!\!\)}}
\DeclareMathOperator\SR{{\rm SR}}

\DeclareMathOperator{\gmetric}{{\mathcal G}}

\numberwithin{equation}{section}
\numberwithin{table}{section}
\def\cc{{\eta}}

\begin{document}

\title[Logarithmic growth of systole]{Logarithmic growth of systole of
arithmetic Riemann surfaces along congruence subgroups}

\def\BIU{Department of Mathematics, Bar Ilan University, Ramat Gan
52900, Israel}

\author[M.~Katz]{Mikhail G. Katz$^{1}$} %
\author[M.~Schaps]{Mary Schaps} %
\author[U.~Vishne]{Uzi Vishne$^2$} %
\address{\BIU} %
\email{\set{katzmik,mschaps,vishne}@math.biu.ac.il} %

\thanks{$^{1}$Supported by the Israel Science Foundation (grants no.\
84/03 and 1294/06)}

\thanks{$^{2}$Supported by the EU research and training network
HPRN-CT-2002-00287, ISF Center of Excellence grant 1405/05, and
BSF grant no.~2004-083}

%\address{}

% \email{}

%\thanks{}



\date{\today}

%\dedicatory{}

\begin{abstract}
We apply a study of orders in quaternion algebras, to the differential
geometry of Riemann surfaces.  The least length of a closed geodesic
on a hyperbolic surface is called its systole, and denoted~$\syspi_1$.
P.~Buser and P.~Sarnak constructed Riemann surfaces~$X$ whose systole
behaves logarithmically in the genus~$g(X)$.  The Fuchsian groups in
their examples are principal congruence subgroups of a fixed
arithmetic group with rational trace field.  We generalize their
construction to principal congruence subgroups of arbitrary arithmetic
surfaces.  The key tool is a new trace estimate valid for an arbitrary
ideal in a quaternion algebra.  We obtain a particularly sharp bound
for a principal congruence tower of Hurwitz surfaces~(PCH), namely the
$4/3$-bound~$\syspi_1(X_{\PCH}) \geq \tfrac{4}{3}\log(g(X_{\PCH}))$.
Similar results are obtained for the systole of hyperbolic
$3$-manifolds, relative to their simplicial volume.
\end{abstract}

\maketitle
\tableofcontents

%-----------------------------------------------------------------

\section{Orders in quaternion algebras and Riemann surfaces}

Arithmetic lattices, besides their own intrinsic interest, have
traditionally provided a rich source of examples in geometry.  One
striking application is the construction of isospectral, non-isometric
hyperbolic surfaces by M.-F. Vigneras \cite{Vi}.  A survey of
arithmeticity as applied in geometry and dynamics may be found in
\cite{Pa}. See \cite{Lub} for an application of congruence subgroups
and the literature on girth in graph theory initiated by
W.~Tutte~\cite{Tu}.  See also \cite{CW} for a recent geometric
application of congruence subgroups.

While the simplest definition of arithmeticity, in analogy
with~$\SL[2](\Z)$, can be presented in terms of~$n$-dimensional
representations by matrices defined over the integers, for many
purposes it is convenient to work with a definition in terms of
quaternion algebras.  The latter is equivalent to the former; \cf
Definition~\ref{14}.  We start by recalling the relevant material on
quaternion algebras.

Let~$a,b\in \overline\Q$, and let
\begin{equation}
\label{11b}
D = \Q[i,j\subjectto i^2=a,\,j^2=b,\,ji=-ij]
\end{equation}
be an (associative) division algebra.  If~$a$ and~$b$ are positive
integers, we have~$D \tensor \R \isom \M[2](\R)$, \cf
\cite[Theorem~5.2.1(i)]{Kato}.  Consider the group~$\Gamma$, by
definition composed of the elements of norm one in~$\Z[i,j] \sub D$.
Then~$\Gamma$ is a co-compact lattice in~$\SL[2](\R)$, see
\cite[Theorem~5.5]{PR}.

P.~Buser and P.~Sarnak \cite[p.~44]{BS} showed that in such case, the
principal congruence subgroups of the Fuchsian group~$\Gamma$ exhibit
near-optimal asymptotic behavior with regard to their systole.  (A
more general construction, but still over~$\Z$, was briefly described
by M.~Gromov \cite[3.C.6]{Gr2}.)  Namely, there is a constant~$c$
independent of~$m$ such that the compact hyperbolic Riemann surfaces
defined as the quotients~$X_m = \dom{\Gamma(m)}{\HC^2}$, satisfy the
bound
$$\syspi_1(X_m) \geq \frac{4}{3} \log \genus(X_m) - c,$$
where~$\HC^2$ is the Poincar\'e upper half plane,~$\Gamma(m)$ are the
principal congruence subgroups of~$\Gamma$,~$\genus(X)$ denotes the
genus of~$X$, and the systole (or girth)~$\syspi_1(X)$ is defined as
follows.
\begin{defn}
The homotopy~$1$-systole, denoted~$\syspi_1(\gmetric)$, of a
Riemannian manifold~$(X,\gmetric)$ is the least length of a
noncontractible loop for the metric~$\gmetric$.
\end{defn}

\begin{rem}
The calculation of \cite{BS} relies upon a lower bound for the
(integer) trace resulting from a congruence relation modulo a rational
prime~$p$.  Such a congruence argument does not go over directly to a
case when the structure constants~$a,b$ of the quaternion
algebra~\eqref{11b} are algebraic integers of a proper extension
of~$\Q$, since the latter are dense in~$\R$. In particular, the
results of~\cite{BS} do not apply to Hurwitz surfaces.  Thus a new
type of trace estimate is needed, see \Tref{boundtr} below.
% \eqref{34} -- ?
\end{rem}

Riemann surfaces with such logarithmic asymptotic behavior of the
systole were exploited by M.~Freedman in his construction of
$(1,2)$-systolic freedom, in the context of quantum computer error
correction~\cite{Fr}.  For additional background on quaternion
algebras and arithmetic Fuchsian groups, see \cite[Appendix to chapter
1]{GGP}, \cite[Section~5.2]{Kato} and~\cite{MR}.  See also the recent
monograph~\cite{SGT} for an overview of systolic problems, as well as
\cite{CK}.

\begin{defn}
\label{14} Let~$F$ denote one of the fields~$\R$ or~$\C$.  An
{\em arithmetic lattice\/}~$G\subset \SL[2](F)$ is a finite co-volume
discrete subgroup, which is commensurable with the group of elements
of norm one in an order of a central quaternion division algebra~$D$
over a number field~$K$ (which has at least one dense embedding
in~$F$).
\end{defn}

\begin{rem}
By a result of K.~Takeuchi~\cite{Ta}, a lattice~$G$ is arithmetic if
and only if the lattice~$G^{(2)}$ (the subgroup generated by the
squares in~$G$) is contained in an order of a division algebra.
\end{rem}

Such a lattice~$\Gamma\subset D$ is cocompact if and only if the
algebra~$D$ as in~\eqref{11b} splits in a single Archimedean place of
the center~$K$ of~$D$ (real or complex depending on~$F$), and remains
a division algebra in all other Archimedean places.

In particular if~$F = \R$ then~$K$ is totally real, and if~$F =
\C$ then~$K$ has one complex place and~$\dimcol{K}{\Q}-2$ real
ones. Denote by~$\algint{K}\subset K$ its ring of algebraic
integers.  Given an ideal~$I\subset \algint{K}$, consider the
associated congruence subgroup~$\Gamma(I) \subset \Gamma$ (see
Definition~\ref{21} for more details).

In the Fuchsian case, we have the following theorem.

\begin{thm}
\label{12b}
Let~$\Gamma$ be an arithmetic cocompact subgroup of\/~$\SL[2](\R)$.
Then for a suitable constant~$c = c(\Gamma)$, the principal congruence
subgroups of~$\Gamma$ satisfy
\[
\syspi_1(X_I) \geq \frac{4}{3} \log \genus(X_I) - c ,
\]
for every ideal~$I \normali \algint{K}$, where~$X_I = \dom {\Gamma(I)}
{\HC^2}$ is the associated hyperbolic Riemann surface.
\end{thm}

The Buser-Sarnak result mentioned above corresponds to the case
where~$D$ is a division algebra defined over~$\Q$, while~$I =
\ideal{m} \normali \Z$.

The proof of \Tref{12b} is given in Section \ref{finish}, where we
provide additional details concerning the constant~$c$, \cf
Theorem~\ref{51}.
% The constant~$c$ is given in the proof ((except for the measure issue))...
%

%\begin{rem}
%It would be helpful to have an explanation in terms of strong
%approximation.
%\end{rem}
%
%% The strong approximation argument is now hidden in Section 4
%% (the mapping modulo the congruence subgroups being onto for
%%  almost all primes), so there is no need to mention it here.

When~$\Gamma(I)$ is torsion free, the area of~$X_I$ is equal to
$4\pi(\genus(X_I)-1)$.  We can therefore rephrase the bound in terms
of the systolic ratio, as follows.

\begin{defn}
The {\em systolic ratio},~$\SR(X, \gmetric)$, of a
metric~$\gmetric$ on an~$n$-manifold~$X$ is
\[
\SR(X, \gmetric)=\frac {\syspi_1(\gmetric)^n} {\vol_n(\gmetric)}.
\]
\end{defn}

\begin{cor}
\label{12c}
% If the ideal~$I\normali \algint{K}$ is large enough in terms of the norm~$\Norm(I)$,
%% - this restriction is only needed when the constant is specified; otherwise there are only finitely many exceptions
%%  (in terms of having torsion).
The hyperbolic surfaces~$X_I = \dom{\Gamma(I)}{\HC^2}$ satisfy the
bound
\[
\SR(X_I) \geq \frac{4}{9\pi} \frac{ \left(\log \genus(X_I) -
c\right)^2} {\genus(X_I)},
\]
where~$c$ only depends on~$\Gamma$.
\end{cor}

Note that an asymptotic {\em upper\/} bound of~$\frac{1}{\pi} \frac
{(\log \genus)^2}{\genus}$ was obtained by the first author in
collaboration with S.~Sabourau in \cite{KS2}, for the systolic ratio
of arbitrary (not necessarily hyperbolic) metrics on a genus~$g$
surface.  The asymptotic multiplicative constant therefore lies in the
interval
\begin{equation}
\label{11} \limsup_{g\to \infty} \frac{g \SR(\Sigma_g)}{\log(g)^2} \in
\left[\tfrac{4}{9\pi} , \tfrac{1} {\pi} \right],
\end{equation}
where~$\Sigma_g$ denotes a surface of genus~$g$.  The question of the
precise asymptotic constant in \eqref{11} remains open. Note that
R.~Brooks and E.~Makover \cite[p.~124]{BM} opine that the Platonic
surfaces of \cite{Br} have systole on the order of~$C \log
g$. However, the text of~\cite{Br} seems to contain no explicit
statement to that effect.  The result may be obtainable by applying
the techniques of~\cite{Br} so as to compare the systole of compact
Platonic surfaces, and the systole of their noncompact prototypes,
namely the finite area surfaces of the congruence subgroups of the
modular group, studied by P.~Schmutz~\cite{Sc}.

The expected value of the systole of a random Riemann surface turns
out to be independent of the genus \cite{MM} (in particular, it does
not increase with the genus), indicating that one does not often come
across surfaces constructed in the present paper.

Another asymptotic problem associated with surfaces is Gromov's
filling area conjecture, when the circle is filled by a surface of an
arbitrary genus~$g$.  The case~$g=1$ was recently
settled~\cite{BCIK1}.

Similarly, in the Kleinian case we have the following.  The simplicial
volume~$\| X \|$, a topological invariant of a manifold~$X$, was
defined in~\cite{gro81}.

\begin{thm}
\label{12d}
Let~$\Gamma$ be an arithmetic cocompact torsion free subgroup
of~$\SL[2](\C)$. Then for a suitable constant~$c = c(\Gamma)$, the
congruence subgroups of~$\Gamma$ satisfy
\[
\syspi_1(X_I) \geq \frac{2}{3} \log \| X_I \| - c ,
\]
for every ideal~$I \normali \algint{K}$, where~$X_I = \dom{\Gamma(I)}
{\HC^3}$, while~$\HC^3$ is the hyperbolic~$3$-space.
\end{thm}

As a consequence, we obtain the following lower bound for the systolic
ratio of the hyperbolic 3-manifolds~$X_I$.  This bound should be
compared to Gromov's similar {\em upper\/} bound, \cf
\cite[6.4.D$''$]{Gr1} (note the missing exponent~$n$ over the log),
\cite[3.C.3(1)]{Gr2}, \cite[p.~269]{Gr3}.

\begin{cor}
\label{16}
Let~$\Gamma$ be an arithmetic cocompact subgroup of~$\SL[2](\C)$.
%
% For~$I$ large enough in terms of the norm~$\Norm(I)$,
%% - Again, this is not needed since the constant is not specified; UV.
%
The~$3$-manifolds~$X_I = \dom{\Gamma(I)}{\HC^3}$ satisfy the bound
\[
\SR(X_I) \geq C_1 \frac{(\log \|X_I\| - c)^3}{\|X_I\|},
\]
where~$c$ and~$C_1$ only depend on~$\Gamma$.  In fact, if~$\Gamma$
is torsion free, then one can take~$C_1 = \frac{8}{27} v_3^{-1}$,
where~$v_3$ is the volume of a regular ideal~$3$-simplex
in~$\HC^3$.
\end{cor}

Theorem~\ref{12d} and Corollary~\ref{16} are proved in
Section~\ref{finish}.  These bounds are shown in \cite{VV} to be exact
if~$K$ has only one Archimedean place.

Asymptotic lower bounds for the systolic ratio in terms of the Betti
number are studied in \cite{BB}. Analogous asymptotic estimates for
the conformal 2-systole of 4-manifolds are studied in \cite{Ka3, Ham}.
See \cite{e7} for a recent study of the~$4$-systole.

In a direction somewhat opposite to ours, hyperbolic~$4$-manifolds of
arbitrarily short~$1$-systole are constructed by I.~Agol \cite{Ag}.

\medskip

Returning the the 2-dimensional case, recall that the order of the
automorphism group of a Riemann surface of genus~$\genus$ cannot
exceed the bound~$84(\genus-1)$. Of particular interest are surfaces
attaining this bound, which are termed Hurwitz surfaces, \cf
\cite{El0, El}.
%
% "The triangle is not a surface since it has boundary.  Its double
% is H/\Gamma, which is a surface with conical singularities
% corresponding to the elliptic elements of~$\Gamma$.  The term
% Hurwitz surface is usually reserved for smooth surfaces.  These
% correspond to subgroups of \Gamma not containing elliptic
% elements. " (MK)
%
Consider the geodesic triangle with angles
$\frac{\pi}{2},\frac{\pi}{3},\frac{\pi}{7}$, which is the least area
triangle capable of tiling the hyperbolic plane. Let~$\Deltah$ be the
group of even products of reflections in the sides of this
triangle. The area of~$\dom{\Deltah}{\HC^2}$, namely~$\pi / 21$, is the
smallest possible for any Fuchsian group.  The Fuchsian group~$N$ of a
Hurwitz surface is a normal torsion free subgroup of~$\Deltah$.  The
automorphism group of the surface is the quotient group~$\Deltah/N$.
The geometry of Hurwitz surfaces was recently studied in~\cite{Vo}.

We specialize to Hurwitz surfaces in Section~\ref{237}.  Following
N.~Elkies \cite{El0}, we choose an order in a suitable quaternion
algebra, as well as a realization of the~$\Z_2$-central
extension~$\tDeltah$ of~$\Deltah$ as a group of~$2\times 2$
matrices.  We then obtain the following sharpening of
Theorem~\ref{12b}.

\begin{thm}\label{special}
For infinitely many congruence subgroups~$\Deltah(I) \normali
\tDeltah$, the Hurwitz surfaces~$X_I = \dom{\Deltah(I)}{\HC^2}$
satisfy the bound
\[
\syspi_1(X_I) \geq \frac{4}{3} \log \genus(X_I).
\]
With the explicit realization of~$\tDeltah$ described in
Section~\ref{237}, the bound holds for all the principal congruence
subgroups.
\end{thm}

Riemannian geometers have long felt that surfaces which are optimal
for the systolic problem should have the highest degree of symmetry;
see, for example, the last paragraph of the introduction to
\cite[p.~250]{HK}.  This sentiment is indeed borne out by our
Theorem~\ref{special} (though in principle there could exist surfaces
with even better asymptotic systolic behavior).

In Section~\ref{B}, we present the key trace estimate which allows us
to generalize the results of Buser and Sarnak to quaternion algebras
over an arbitrary number field.  In Section~\ref{three}, we prove the
trace estimate.  Section~\ref{sec:cong} contains a detailed study of
congruence subgroups. We comment on torsion elements in
Section~\ref{five}.  Section~\ref{finish} contains the proofs of the
main theorems.  Section~\ref{237} focuses on the Hurwitz case.

\section{Trace estimate}
\label{B}

Fix~$F = \R$ or~$F = \C$, and let~$K$ be a number field of
dimension~$d$ over~$\Q$, as in Definition~\ref{14}. Namely,~$K$ is a
totally real field in the former case, or a field with exactly one
complex Archimedean place, in the latter.  We view~$K$ as a subfield
of~$F$, where the other real embeddings are denoted by~$\s\co K \hra
\R$ for~$\s\neq 1$. Let~$D$ be a quaternion division algebra with
center~$K$.  We assume that~$D$ is split by the distinguished
embedding~$K \sub F$, and remains a division algebra under any other
embedding.

Recall that an {\em order\/} $\CQ$ of~$\quat$ is a subring (with
unit), which is a finite module over $\Z$, and such that its
(central) ring of fractions is~$D$. Every maximal order
contains~$\algint{K}$, the ring of algebraic integers in $K$
(since $\algint{K}\cdot \CQ$ is an order), so in the sequel we
will only deal with orders containing~$\algint{K}$.

Let~$D^1 \sub \mul{D}$ denote the group of elements of norm~$1$
in~$D$. Similarly, for an order~$\CQ \subset D$,
let~$\SLform{\CQ}$ denote the group of elements of norm one
in~$\CQ$.  By assumption the inclusion~$K \sub F$ splits~$D$, so
we have the natural inclusion
\begin{equation}
\label{21b}
\CQ \sub D \sub D \tensor[K]F \isom \M[2](F).
\end{equation}
Thus~$\SLform{\CQ}$ is a subgroup of~$(D\tensor[] F)^1 =
\SL[2](F)$.

Since~$\CQ$ is an order over~$\algint{K}$, in particular it
contains~$\algint{K}$ (as its center).  An ideal~$I$ of the center
defines an ideal~$\IQ$ of~$\CQ$, yielding a finite quotient
ring~$\CQ/\IQ$. The {\it principal congruence subgroup} of~$\CQ^1$
with respect to an ideal~$I \normali \algint{K}$ is by definition the
kernel of the homomorphism~$\CQ^1 \ra \mul{(\CQ/\IQ)}$ induced by the
natural projection~$\CQ \ra \CQ/\IQ$.  A {\it congruence subgroup} is
any subgroup of~$\CQ^1$ containing a principal congruence subgroup.

\begin{defn}
\label{21} The kernel of~$\CQ^1 \ra \mul{(\CQ/\IQ)}$ induced by
the natural projection~$\CQ \ra \CQ/\IQ$, will be denoted
by~$\Gamma(I)= \CQ^1(I)$, where
\begin{equation}
\label{22} \CQ^1(I) = \ker \left( \CQ^1 \ra \mul{(\CQ/\IQ)}
\right) .
\end{equation}
\end{defn}
Since~$\quat \tensor[\s] \R$ is, by hypothesis, a division algebra
for~$\s \neq 1$, the groups~$(\quat \tensor[\s] \R)^1$ are
compact. Moreover~$\algint{K}$ is discrete in the product~$F \times
\prod_{\s \neq 1} \s(\R)$, making~$\CQ^1$ discrete in the product
\[
(D\tensor[] F)^1 \times \prod (\quat\tensor[\s]\R)^1,
\]
and therefore discrete in~$(D\tensor[] F)^1 = \SL[2](F)$. In fact
$\CQ^1$ is cocompact there, by \cite[Theorem~5.5]{PR}.

Since all orders of~$D$ are commensurable, the groups of elements of
norm~$1$ in orders are commensurable to each other, \ie the
intersection of two such groups is of finite index in each of the
groups~\cite[p.~56]{MR}.  Our computations are based on a specific
order, arising naturally from the presentation of~$\quat$.

A quaternion algebra can always be presented in the form %
\begin{equation}\label{Ddef}
\quat = (a,b)_K = K[i,j \subjectto i^2 = a,\, j^2 = b,\, ji =
-ij]%
\end{equation}
for suitable elements~$a,b \in K$. If~$\s \co K \hra \R$ is an
embedding, then with this presentation, the algebra~$\quat\tensor[\s]
\R$ is a division algebra if and only if we have both~$\s(a)<0$
and~$\s(b) < 0$. Another convenient feature is that the matrix trace
in the embedding~$D \sub \M[2](F)$ is equal to the reduced
trace~$\Tr[D]$, which for non-central elements is the negative of the
linear coefficient in the minimal polynomial.  This can be read off
the presentation of an element:
$$\Tr[D](x_0 + x_1 i + x_2 j + x_3 ij) = 2x_0$$
for any~$x_0,x_1,x_2,x_3 \in K$. %
\begin{lemma}
\label{2.1} The defining constants~$a$ and~$b$ of the algebra~$D$
can be taken to be algebraic
integers in~$K$.
% Moreover~$\s(a),\s(b) < 0$ for any~$\s \neq 1$.
\end{lemma}

\begin{proof}
It is clear from the presentation~\eqref{Ddef} that the isomorphism
class of~$D$ depends only the class of~$a$ and~$b$ in~$\md[2]{K}$.
\end{proof}

Let~$\algint{K}$ be the ring of algebraic integers in~$K$,  %
and fix the order %
\begin{equation}
\label{standord}
\O = \algint{K} \oplus \algint{K} i \oplus \algint{K} j \oplus
\algint{K} ij.
\end{equation}

Since all elements of an order are algebraic integers,~$\CQ$ is
contained in~$\frac{1}{\kappa} \O$ for a suitable~$\kappa \in
\algint{K}$ (in fact one can take~$\kappa \divides 2ab$, \cf
Lemma~\ref{2.1}).
Denote by~$\Norm(k)$ the number field norm of~$k\in K$ along the
extension~$K/\Q$. Similarly, we denote by~$\Norm(I)$ the norm of an
ideal~$I \normali \algint{K}$, namely the cardinality of the quotient
ring~$\algint{K}/I$. The two norms coincide for principal ideals, so
in particular, we have~$\Norm(m) = m^d$ for~$m \in \Z$.

\begin{thm}
\label{boundtr}
Let~$I \normali \algint{K}$ be an ideal.
% , and~$\CQ\sub D$ and~$\kappa$ as before.
If~$F = \R$ then for every~$x \neq \pm 1$ in~$\Lat(I)$, we have
the following estimates:
\begin{equation}
\abs{\Tr[D](x)} > \frac{1}{2^{d-2}\Norm(\ideal{2} + \kappa I)}
\Norm(I)^2 - 2
\end{equation}
and therefore
\begin{equation}
\abs{\Tr[D](x)} >  \frac{1}{2^{2d-2}} \Norm(I)^{2} - 2.
\end{equation}
If~$F = \C$, then
\begin{eqnarray*}
\abs{\Tr[D](x)}
    & > &
        \frac{1}{2^{d/2 - 2}
            \Norm(\ideal{2}+\kappa I)^{1/2}}{\Norm(I)} - 2 \\
    & \geq &
        \frac{1}{2^{d-2}}\abs{\Norm(I)} - 2.
\end{eqnarray*}
\end{thm}
The estimates are proved in the next section.

\section{Proof of trace estimate}
\label{three}

The first step towards a proof of Theorem~\ref{boundtr} is a bound
on~$\s(x_0)$ for~$\s\neq 1$, where~$x \in D$ is an arbitrary element
of norm~$1$.

\begin{Prop}\label{boundsx}
Let~$x = x_0 + x_1 i + x_2 j + x_3 ij$, where~$x \neq \pm 1$, be any
element of norm one in~$\quat$.  Then~$\abs{\s(x_0)} < 1$ for every
non-trivial embedding~$\s \co K \hra \R$. Writing~$x_0 = 1 + y_0$, we
obtain
\begin{equation}
\label{ineq}
-2 < \s(y_0) < 0.
\end{equation}
\end{Prop}
\begin{proof}
The norm condition is that %
\begin{equation}
\label{xx}%
x_0^2 - a x_1^2 - b x_2^2 + ab x_3^2 = 1
\end{equation}
in~$K$, where~$a,b$ are the structure constants from~\eqref{Ddef}.
Applying any~$\s \neq 1$, we obtain
\begin{equation}
\begin{aligned}
\s(x_0)^2 & \leq \s(x_0)^2 + \s(- a) \s(x_1)^2 + \s(-b) \s(x_2)^2
+ \s(ab) \s(x_3)^2
\label{33}
\\&= 1,
\end{aligned}
\end{equation}
since~$\s(a)<0$ and~$\s(b) < 0$ by assumption.  In particular,
equality in~\eqref{33} implies that~$x_1 = x_2 = x_3 = 0$.
Writing~$x_0 = 1 + y_0$, we obtain the inequality~$\abs{1+\s(y_0)} <
1$, proving \eqref{ineq}.
\end{proof}

Recall that for any ideal~$J \normali \algint{K}$, the fractional
ideal
\[
J^{-1} = \set{u \in K \suchthat u J \sub \algint{K}}
\]
is the inverse of~$J$ in the group of fractional ideals of~$K$.

Recall that the symplectic involution on~$D$ is, by definition,
the unique involution under which only central elements are
symmetric. It is often called the `standard' involution.
\begin{lemma}
\label{invpres}
The symplectic involution preserves any order in~$D$.
\end{lemma}
\begin{proof}
Let~$w\mapsto w^*$ denote the involution. The reduced trace and norm
may be defined by~$\Tr[D](w) = w+ w^*$ and~$\Norm[D](w) = ww^*$, which
implies the characteristic equation
\[
w^2 - \Tr[D](w) w + \Norm[D](w) = 0
\]
for every~$w \in D$. Moreover if~$\CQ$ is an order and~$w \in \CQ$, it
follows that~$\Tr[D](w),\Norm[D](w) \in \algint{K}$, as~$\algint{K}[w]
\sub \CQ$ is a finite module.  In particular, we have~$w^* = \Tr[D](w)
- w \in \CQ$ for every~$w \in \CQ$.
\end{proof}

\begin{lemma}\label{32}
Let~$I \normali \algint{K}$. If~$z \in \IQ$ then~$\Tr[D](z)\in
I$ and~$\Norm[D](z) \in I^2$.
\end{lemma}
\begin{proof}
Let~$\alpha_1,\dots,\alpha_t \in I$ be generators of~$I$ as an
$\algint{K}$-module (in fact we may assume~$t \leq 2$,
as~$\algint{K}$ is a Dedekind domain \cite[Section~VII.10]{Cohn}).

Let~$z = \sum \alpha_r w_r$,~$w_r \in \CQ$, be an arbitrary element
of~$\IQ$.  Then we have~$w_r w_s^* \in \CQ$ by Lemma~\ref{invpres}, so
that
\[
\Tr[D](z) = \sum\alpha_r \Tr[D](w_r) \in I
\]
since~$\Tr[D](w_r) \in \algint{K}$, and
\begin{eqnarray*}
\Norm[D](z) = z z^* & = & \sum_{r,s} \alpha_r \alpha_s w_r w_s^* \\
    & = & \sum_{r} \alpha_r^2 w_r w_r^* + \sum_{r < s} \alpha_r
    \alpha_s (w_r w_s^* + w_s w_r^*) \\
    & = & \sum_{r} \alpha_r^2 \Norm[D](w_r) + \sum_{r < s}
    \alpha_r
    \alpha_s \Tr[D](w_r w_s^*) \in I^2
\end{eqnarray*}
since~$\Norm[D](w_r), \Tr[D](w_r w_s^*) \in \algint{K}$.
\end{proof}

\begin{lemma}
Let~$x = x_0 + x_1 i + x_2 j + x_3 ij \in \CQ^1(I)$, where~$\CQ \sub
\frac{1}{\kappa}\O$, and~$\O$ is the standard order
of~\eqref{standord}.  Then
\begin{equation}
\label{whereisy} %
x_0-1\in  (\ideal{2} + \kappa I)^{-1} I^2.
%%% Used to be:~$\frac{1}{\kappa}(\ideal{2\kappa}+I) I^2$.  %
\end{equation} %
\end{lemma}
\begin{proof}
By assumption~$x -1 \in \IQ \sub \frac{1}{\kappa}I\O$. In
particular, the coefficient~$y_0 = x_0-1 \in \frac{1}{\kappa}I$.

Substituting~$x_0 = 1 + y_0$ in \eqref{xx}, we obtain
\begin{equation}\label{xx2}
2 y_0 = - (y_0^2 - a x_1^2 - b x_2^2 + ab x_3^2) = - \Norm[D](x-1)
\in I^2,
\end{equation}
by the previous lemma.  Using the decomposition of fractional
ideals in the Dedekind domain~$\algint{K}$, which implies
that~$(J_1+J_2)(J_1\cap J_2) = J_1J_2$, we obtain
\begin{equation}
\label{34}
y_0 \in \frac{1}{\kappa} I \cap \frac{1}{2} I^2 = \frac{1}{2\kappa} (2
\algint{K} \cap \kappa I) I = (2\algint{K}+\kappa I)^{-1} I^2,
\end{equation}
proving the lemma.
\end{proof}

\begin{proof}[Proof of Theorem~\ref{boundtr}]
We write~$x = x_0 + x_1 i + x_2 j + x_3 ij$, and we set~$y_0 = x_0 -
1$.  By the lemma, we have~$\Norm(y_0) \in
\frac{\Norm(I)^2}{\Norm(\ideal{2} + \kappa I)}\Z$, but we also
have~$y_0 \neq 0$ (\Pref{boundsx}). When~$F = \R$, this shows
\begin{eqnarray*}
\frac{\Norm(I)^2}{\Norm(\ideal{2}+\kappa I)} \leq \abs{\Norm(y_0)}
& = & \abs{y_0} \prod_{\s \neq 1} \abs{\s(y_0)}
     <  2^{d-1}\abs{y_0},
\end{eqnarray*}
by inequality \eq{ineq}. The inclusion~$\ideal{2} \sub
\ideal{2}+\kappa I$ implies
\[
 \Norm(\ideal{2}+ \kappa I) \leq \Norm(\ideal{2}) = \Norm(2)
 = 2^d.
\]
Plugging the norm bound on~$y_0$ back in~$x_0 = y_0+1$, we obtain
\begin{eqnarray*}
\abs{x_0}  & =  & \abs{1+y_0} \nonumber \\
     & \geq &  \abs{y_0} - 1 \label{realbound} \\
     & > & \frac{1}{2^{d-1}
        \Norm(\ideal{2} + \kappa I)} \Norm(I)^{2} -
        1,\nonumber \\
     & \geq & \frac{1}{2^{2d-1}} \Norm(I)^{2} - 1,\nonumber
\end{eqnarray*}
from which the case~$F = \R$ of \Tref{boundtr} follows immediately
in view of the fact that~$\Tr[\quat](x) = 2x_0$.

In the case~$F = \C$, we repeat the same argument, noting that
$\abs{\Norm(y_0)} = \abs{y_0}^2 \prod_{\s \neq 1} \abs{\s(y_0)}$
so that~$ \frac{\abs{\Norm(I)}^2}{\Norm(\ideal{2} + \kappa I)} <
2^{d-2}\abs{y_0}^2.~$
\end{proof}

\section{Congruence subgroups of~$\Lat$ in the prime power case}
\label{sec:cong}

Some of the calculations in this section are related to A.~Borel's
volume formula \cite{Bo}, see also \cite{Jo, Be}.

Let~$F$ be~$\R$ or~$\C$ as above, and let~$\Gamma = \Lat$ be the group
defined by an order~$\CQ$ in a division quaternion algebra~$\quat$
over a number field~$K$, as in Sections \ref{B} and~\ref{three}.

Let~$I \normali \algint{K}$ be an arbitrary ideal.  In order to
estimate the topological invariants of~$\dom{\Gamma(I)}{\HC^2}$ or
$\dom{\Gamma(I)}{\HC^3}$ (namely, the genus and simplicial volume,
respectively), we need to bound the index~$\dimcol{\Gamma}{\Gamma(I)}$
in terms of~$I$.  We will construct below a norm map
\begin{equation}
\label{41} \nu\co\mul{(\CQ/\IQ)} \ra \mul{(\algint{K}/I)}.
\end{equation}
Let~$(\CQ/\IQ)^1= \nu^{-1}(1)$, where~$1$ denotes the
multiplicative neutral element in the finite ring~$\algint{K}/I$.

\begin{lemma}
We have the bound
\begin{equation}
\label{42} \dimcol{\Gamma}{\Gamma(I)} \leq (\CQ/\IQ)^{1}.
\end{equation}
\end{lemma}

\begin{proof}
The reduced norm~$N_\quat$ on~$D$ may be defined by setting~$N(x) =
xx^*$, where~$x\mapsto x^*$ is the unique symplectic involution
of~$D$.  Note that every element with zero trace is anti-symmetric
under the involution.  After tensoring with a splitting field, the
involution becomes the familiar adjoint map, which exchanges the two
diagonal elements and multiplies the off diagonal elements by~$-1$.
Note that orders are closed under involution, as~$x^{*} = \tr(x)-x$.
% ($x^*+x$ is the trace)

Note that the symplectic involution descends to the
quotient~$\CQ/\IQ$.  Indeed, the involution acts trivially
on~$I\normali \algint{K}$ and preserves the order~$\CQ \subset D$ by
Lemma~\ref{invpres}.  Hence an ideal of the form~$\IQ$ is also closed
under the involution.  Thus we may define the involution on a
coset~$(x+\IQ)$ by setting~$(x+\IQ)^* = x^* + \IQ$, which is
independent of the representative. The involution can then be used to
define a norm map~$\nu$ of \eqref{41}. Let~$\pi \co \CQ \ra \CQ/\IQ$
and~$\pi_0 \co \algint{K} \ra \algint{K}/I$ denote the natural
projections, then~$\pi_0 \Norm(x) = \pi_0(xx^*) = \nu(\pi(x))$, and in
particular~$\pi(x) \in (\CQ/\IQ)^1$ for every~$x \in \CQ^1 = \Gamma$.

By definition,~$\Ker(\pi) \cap \Gamma = \Gamma(I) = \CQ^1(I)$, and
so
$$\Gamma/\Gamma(I) \isom \pi(\CQ^1) \leq (\CQ/\IQ)^1,$$
proving the claim.
\end{proof}

Since~$K$ is an algebraic number field,~$\algint{K}$ is a Dedekind
domain, so that every ideal factors as an intersection of powers
of prime ideals. Suppose~$I = \fp^t$ where~$\fp \normali
\algint{K}$ is a prime ideal. To find an upper bound for the
right-hand side of \eqref{42}, we pass to the local situation, via
standard techniques in algebraic number theory, as follows.

Let~$K_\fp$ be the completion of~$K$ with respect to the~$\fp$-adic
valuation, and~$O_{\fp}$ the valuation ring.  We consider
an~$O_{\fp}$-order in~$D_{\fp} = D \tensor[K] K_{\fp}$, defined by
setting~$\CQ_{\fp} = \CQ \tensor[{\algint{K}}] O_{\fp}$.

\begin{lemma}
We have an isomorphism~$\CQ/\IQ = \CQ_{\fp} / \fp^t \CQ_{\fp}$.
\end{lemma}

\begin{proof}
Let~$S = \algint{K}\minusset \fp$, the complement of~$\fp$
in~$\algint{K}$. Since~$\fp$ is a maximal ideal, localization (in
which elements of~$S$ are forced into being invertible) gives
\[
\CQ / \fp^t \CQ \isom S^{-1}\CQ / \fp^t S^{-1} \CQ,
\]
which is isomorphic to~$\CQ_{\fp} / \fp^t \CQ_{\fp}$
\cite[Exercise~5.7]{Reiner}.
\end{proof}

It follows that
\begin{equation}\label{bound1}
\dimcol{\Gamma}{\Gamma(\fp^t)} \leq
\card{(\CQ_\fp/\fp^t\CQ_\fp)^{1}}.
\end{equation}

In one common situation, we can compute
$\card{(\CQ_\fp/\fp^t\CQ_\fp)^{1}}$ explicitly. To put things in
perspective, it is worth noting that for all but finitely many
primes~$\fp$, the algebra~$D_{\fp}$ is a matrix algebra,
and~$\CQ_{\fp}$ is maximal, \ie not contained in a larger order.
Let~$q = \Norm(\fp) = \card{\algint{K}/\fp}$.

\begin{lemma}\label{maxim}
Suppose~$\CQ_{\fp}$ is a maximal order of~$D_{\fp}$.  If~$D_{\fp}$ is
a division algebra, then~$\CQ/\fp^t\CQ$ is a local (non-commutative)
ring with residue field of order~$q^2$ and radical whose nilpotency
index is~$2t$; in such case, we have~$\card{(\CQ/\fp^t\CQ)^1} =
q^{3t}(1+\frac{1}{q})$. Otherwise, we have~$\CQ/\fp^t\CQ \isom
\M[2](\algint{\fp}/\fp^t)$, and~$\card{(\CQ/\fp^t\CQ)^1} =
q^{3t}(1-\frac{1}{q^2})$.
\end{lemma}

%% Note: the strucuture can be described in great details, and not only for~$\deg(D)=2$.

\begin{proof}
Let~$\varpi \in O_{\fp}$ be a uniformizer (\ie generator of the unique
maximal ideal).  Central simple algebras over local fields are cyclic
\cite[Chapter~17]{Pierce}.  Therefore, there is an unramified maximal
subfield~$L$, satisfying
\begin{equation}
K_{\fp} \sub L \sub D_{\fp},
\end{equation}
and an element~$z \in D_{\fp}$, such that~$z \ell z^{-1} = \s(\ell)$
for~$\ell \in L$ (where~$\s$ is the non-trivial automorphism
of~$L/K_{\fp}$), and moreover~$D_{\fp} = L[z]$. Furthermore
if~$D_{\fp}$ is a division algebra then we may assume~$z^2 = \varpi$.
If~$D_{\fp}$ splits, then being a matrix algebra it can be presented
in a similar manner, with~$z^2 = 1$.

In both cases,~$D_\fp$ has, up to conjugation, a unique maximal
order \cite[Section~17]{Reiner}, namely~$O_L[z]$ where~$O_L$ is
the valuation ring of~$L$.
% Since every order is contained in a
% maximal order, we may assume~$\CQ_\fp \sub O_L[z]$. Thus
%~$\CQ_{\fp} / \fp^t \CQ_{\fp} = \CQ_{\fp} /\varpi^t \CQ_{\fp}$ is
% isomorphic to a subalgebra of~$O_L[z]/\varpi^t O_L[z]$.
We may therefore assume~$\CQ_{\fp} = O_L[z]$.

To finish the proof, we note that~$O_L[z]/\varpi^t O_L[z]$ is a
local ring, so its invertible elements are those invertible modulo
the maximal ideal. First assume~$D_{\fp}$ is split, so that~$z^2 =
1$. Then
\[
O_L[z]/\varpi^t O_L[z] = (O_L/\varpi^t O_L)[\bar{z} \,|\,
\bar{z}^2=1]
\]
is isomorphic to~$\M[2](O_\fp/\varpi^t O_\fp)$, noting
that~$O_L/\varpi^t O_L$ is a finite local ring, which modulo the
maximal ideal is the residue field~$O_\fp/\varpi O_\fp = \bar{K}$, of
order~$q$.  The group of elements of norm~$1$
is~$\SL[2](O_{\fp}/\varpi^t O_{\fp})$, of order~$(q^3-q)q^{3(t-1)}$.

Finally, assume~$D_{\fp}$ does not split, so~$z^2 = \varpi$. Then
the ring
\[
O_L[z]/\varpi^t O_L[z] = O_L[z]/z^{2t} O_L[z]
\]
is a local (non-commutative) ring, which modulo the maximal ideal is
\[
O_L[z]/z O_L[z] \isom \bar{L},
\]
the residue field of~$L$, of order~$q^2$.  Thus there
are~$(q^2-1)q^{2(2t-1)}$ invertible elements.  The norm map
\[
\mul{(O_L[z]/z^{2t} O_L[z])} \ra \mul{(O_\fp/\varpi^t)}
\]
is onto (since the norm~$\mul{O_L} \ra \mul{O_\fp}$ is), and
since~$\card{\mul{(O_\fp/\varpi^t)}} = (q-1)q^{t-1}$, there are
precisely
\[
\frac{(q^2-1)q^{2(2t-1)}}{(q-1)q^{t-1}} = (q+1)q^{3t-1}
\]
elements of norm~$1$.
%
%% There is a messier, more direct way to count elements of norm~$1$;
%% for that we need to write a general element as~$\sum_{\alpha_i z^i}$
%% (with~$\alpha_i$ from a list of~$\bar{L}$ to~$O_L[z]/z^{2t}$) and
%% compute the norm. There are~$q+1$ options for~$\alpha_0$, then
%%~$q^2$ for~$\alpha_1$. From that point,~$q$ options to~$\alpha_2$,
%%~$q^2$ to~$\alpha_3$;~$q$ to~$\alpha_4$,~$q^2$ to~$\alpha_5$, etc.
\end{proof}

Our next goal is to bound~$(\CQ_\fp/\fp^t\CQ_\fp)^{1}$ in the general
case. Consider the exact sequence induced by the norm map, namely
\begin{equation}\label{Nmap} (\CQ_\fp/\fp^t\CQ_\fp)^{1}
\hra (\CQ_\fp/\fp^t\CQ_\fp)^{\times} \stackrel{N}{\longrightarrow}
(\algint{K}/\fp^t)^{\times},\end{equation} which shows that
\begin{equation}\label{bound2}
\card{(\CQ_\fp/\fp^t\CQ_\fp)^{1}} \leq
\frac{\card{(\CQ_\fp/\fp^t\CQ_\fp)^{\times}}}
{\card{\opname{Im}(N)}}.\end{equation} Here,~$\opname{Im}(N)$ stands
for the image of the norm map in \eq{Nmap}. We will first treat the
numerator, then the denominator.

Since~$\fp\CQ_{\fp}/\fp^t\CQ_{\fp}$ is a nilpotent ideal of
$\CQ_\fp/\fp^t\CQ_\fp$, the invertible elements in this finite
ring are those invertible modulo~$\fp\CQ_{\fp}/\fp^t\CQ_{\fp}$;
namely, in the quotient ring~$\CQ_{\fp}/\fp\CQ_{\fp}$. Let~$u =
\card{(\CQ_\fp/\fp\CQ_\fp)^{\times}}$, then
\[
\card{(\CQ_\fp/\fp^t\CQ_\fp)^{\times}} = \card{\fp\CQ_{\fp}/
\fp^t\CQ_{\fp}} u \leq q^{4(t-1)}u.
\]
By Wedderburn's decomposition theorem, we can decompose
the~$\algint{K}/\fp$-algebra
$\CQ_{\fp}/\fp\CQ_{\fp}$ as~$T \oplus J$, where~$T$ is a semisimple
algebra and~$J$ is the radical. Again, the invertible elements are
those invertible modulo~$J$, namely in~$T$.

\begin{prop}
Let~$\F_{q}$ and~$\F_{q^2}$ denote the finite fields of orders~$q$
and~$q^2$, respectively. Then~$T$ is isomorphic to one of the
following six algebras:
\[
\F_q,\ \F_{q^2},\ \M[2](\F_q),\ \F_q\times \F_q, \ \F_q \times
\F_{q^2},\ \mbox{or} \quad \F_{q^2}\times \F_{q^2}.
\]
Furthermore, only the first three options are possible if~$D_{\fp}$ is
a division algebra.
\end{prop}
\begin{proof}
By construction,~$T$ is a semisimple algebra over~$\F_q$. The
classification follows from two facts: every element satisfies a
quadratic equation; and the maximal number of mutually orthogonal
idempotents is no larger than the corresponding number for~$\CQ_{\fp}$
(as Hensel's lemma allows lifting idempotents from~$\CQ/\fp\CQ$
to~$\CQ$).
\end{proof}

By inspection of the various cases, the number~$u$ of invertible
elements satisfies
\begin{equation}
\label{bound3} u \leq
\begin{cases}q(q-1)(q^2-1) & \mbox{if~$D_{\fp}$ is a division
algebra}, \\q^2(q-1)^2 & \mbox{otherwise}.\end{cases}
\end{equation}

It remains to give a lower bound for the size of~$\opname{Im}(N)$.
Since~$\algint{\fp} \sub \CQ$, every square in
$(\algint{\fp}/\fp^t)^{\times}$ is a norm.

\begin{lemma}
\label{44} If~$\fp$ is non-diadic (namely~$\fp$ is prime to~$2$,
and~$q$ is odd), then~$\card{{(\algint{K}/\fp^t)^{\times}}^2} =
\frac{q-1}{2}q^{t-1}$.

On the other hand if~$\fp$ is diadic and~$2 \algint{\fp} =
\fp^{e}$ for~$e\geq 1$, then
$\card{{(\algint{K}/\fp^t)^{\times}}^2} \geq \frac{1}{2}q^{t-e}$
(with equality if~$t\geq 2e+1$).
\end{lemma} % -- ??

\begin{proof}
Note that an element ~$a \in \algint{\fp}^{\times}$ is a square if and
only if~$a$ is a square modulo~$4\fp$ \cite[Chapter~63]{OM}.
Furthermore, we have~$y^2 \equiv x^2 \pmod{4\fp}$ if and only if~$y
\equiv \pm x \pmod{2\fp}$, proving the lemma.
\end{proof}

Combining Equations~\eq{bound2} and \eq{bound3} with
Lemma~\ref{44}, we have the following bounds for non-maximal
orders:
\begin{equation}
\label{bound3bis}
{}\frac{(\CQ/\fp^t\CQ)^1}{q^{3t}} \leq
\begin{cases}%
2(1-q^{-2}) & \mbox{if~$D_{\fp}$ is a division algebra,~$\fp$ non-diadic}, %
\\ %
2(1-q^{-1})(1-q^{-2})q^{e} & \mbox{if~$D_{\fp}$ is a division algebra,~$\fp$ diadic}, %
\\ %
2(1-q^{-1}) & \mbox{if~$D_{\fp}$ is a matrix algebra,~$\fp$ non-diadic}, %
\\ %
2(1-q^{-1})^2q^{e} & \mbox{if~$D_{\fp}$ is a matrix algebra,~$\fp$ diadic}, %
\end{cases}
\end{equation}
where \begin{equation}\label{defe} 2\algint{\fp} =
\fp^e\end{equation} in the diadic cases.
% To illustrate the difficulty, one may consider the order~$\CQ_\fp
% = \algint{\fp}+\fp^{t}M$, where~$M$ is the maximal order over~$\algint{\fp}$.

\bigskip

Given an order~$\CQ$ of a quaternion algebra~$D$ over~$K$, let~$T_1$
denote the set of finite primes~$\fp$ for which~$D_{\fp}$ is a
division algebra, and let~$T_2$ denote the set of finite primes for
which~$\CQ_{\fp}$ is non-maximal. It is well known that~$T_1$
and~$T_2$ are finite.

We denote
\begin{equation}\label{lamdef}
\lam_{D,\CQ} = \prod_{\fp \in T_1\minusset T_2}
{\(1+\frac{1}{\Norm(\fp)}\)} \cdot \prod_{\fp \in T_2} 2 \cdot
\prod_{\fp \in T_2, \fp\divides 2} \Norm(\fp)^{e(\fp)},
\end{equation}
where for a diadic prime,~$e(\fp)$ denotes the ramification index
of~$2$, as defined in \eqref{defe}. The third product is bounded from
above by
\[
\prod_{\fp\divides 2} \Norm(\fp)^{e(\fp)} = \Norm(2) = 2^d .
\]
Factoring an arbitrary ideal~$I \normali \algint{K}$ into prime
factors~$I = \prod {\fp_i^{t_i}}$, we have by the Chinese remainder
theorem that~$\CQ/\IQ \isom \prod \CQ/\fp_i^{t_i}\CQ$, where the
projection onto each component preserves the norm (as the norm can be
defined in terms of the involution).  The results of Lemma~\ref{maxim}
and \eqref{lamdef} imply the bound
\[
\card{(\CQ/\IQ)^1} \leq \lam_{D,\CQ} \prod \Norm(\fp_i^{t_i})^3 =
\lam_{D,\CQ} \Norm(I)^3.
\]
Together with Equation~\eq{bound1}, we obtain the following corollary.

\begin{cor}
\label{3.2} We have the bound~$\dimcol{\Gamma}{\Gamma(I)} \leq
\card{(\CQ/\IQ)^1} \leq \lam_{D,\CQ} \Norm(I)^{3}$. In fact, the
products in \eqref{lamdef} can be taken over the primes~$\fp \divides
I$; in particular if~$\CQ$ is maximal and no ramification primes
of~$D$ divide~$I$, then~$\dimcol{\Gamma}{\Gamma(I)} < \Norm(I)^{3}$.
\end{cor}

\section{Torsion elements}
\label{five}

A comment about torsion elements in~$\Lat$ is in order. If~$x \in
\quat$ is a root of unity, then~$K[x]$ is a quadratic field
extension of~$K$, and with~$K$ fixed, there are of course only
finitely many roots of unity with this property. Moreover if~$x
\in \CQ^1$ then~$\Norm(x) = 1$ forces~$x+x^{-1} \in K$. Now, if~$x
\in \CQ^1(I)$ for an ideal~$I \normali \algint{K}$, then we have
\[
0 = x^2 - (x+x^{-1})x +1 \equiv 2-(x+x^{-1}) = -x^{-1}(x-1)^2
\pmod{I},
\]
so~$I \supseteq \ideal{x^{-1}(x-1)^2}$. These are the ideals to be
avoided if we want~$\CQ^1(I)$ to be torsion free:
\begin{cor}\label{torf}
If~$\Lat(I)$ is not torsion free, then~$I$ divides an ideal of the
form
\[
\ideal{x^{-1}(x-1)^2} = \ideal{x+x^{-1}-2},
\]
when~$x$ a root of unity for which~$K[x]/K$ is a quadratic extension.
Moreover if~$I$ contains a principal ideal~$I_0$, then~$I_0^2 \divides
\ideal{x^{-1}(x-1)^2}$.
\end{cor}
\begin{proof}
Only the final statement was not proved. Suppose~$x \in \Lat(I)$
is a root of unity; then~$K[x]/K$ is a quadratic extension, and
$x^2 - (x+x^{-1})x + 1 =0$. Thus~$\Tr[D](x) = x+x^{-1} \in
\algint{K}$, and~$\Norm(x) = 1$.

Suppose~$I_0 \sub I$ and~$I_0 = \ideal{i}$ is principal.  Write~$x =
1+i a$ for~$a \in \CQ$.  Then
\[
x+x^{-1} = \Tr[D](x) = 2+i \Tr[D](a),
\]
so that~$\Tr[D](a) =
\frac{x+x^{-1}-2}{i}$, and
\[
1 = \Norm(x) = (1+ia)(1+ia^*) = 1+i(a+a^{*}) + i^2 aa^{*},
\]
where~$a^*$ is the quaternion conjugate.  Therefore,
\[
\Norm(a) = aa^* = - \frac{1}{i}(a+a^*) = \frac{2-x+x^{-1}}{i^2}.
\]
But as an element of an order, the norm of~$a$ is an algebraic
integer, implying that~$i^2$ divides~$(x+x^{-1})-2$.
\end{proof}

\section{Proof of the main theorems}
\label{finish}

In this section we prove Theorems \ref{12b} and \ref{12d}.  First,
we recall the relation between the length of closed geodesics and
traces. Let~$\Gamma \leq \SL[2](F)$ be an arbitrary discrete
subgroup, and set~$X = \dom{\Gamma}{\HC^2}$ (or~$X =
\dom{\Gamma}{\HC^3}$ when~$F=\C$). Let~$x \in \Gamma$ be a
semisimple element.  Then~$x$ is conjugate (in~$\SL[2](F)$) to a
matrix of the form~$\smat{\lam}{0}{0}{\lam^{-1}}$, and
\[
\abs{\lam} = e^{\ell_x/2},
\]
where~$\ell_x > 0$ is the translation length of~$x$, which is the
length of the closed geodesic corresponding to~$x$ on the manifold
$X$.  Note that if~$F =
\R$, then in fact~$\lam = e^{\ell_x/2}$. In the Kleinian case~$F=\C$,
there is a rotation of the loxodromic element, as well. In either
case, we have
$$
\abs{\Tr[{\M[2](F)}](x)} = \abs{\lam+\lam^{-1}} \leq \abs{\lam} +
\abs{\lam}^{-1} < \abs{\lam} + 1,
$$
and
\begin{equation}\label{boundell}
\ell_x > 2 \log(\abs{\Tr[{\M[2](F)}](x)} - 1).
\end{equation}

Because of the undetermined constant in theorems \ref{12b} and
\ref{12d}, it is enough to treat the co-compact lattice~$\Gamma$
of~$\SL[2](F)$ up to commensurability. However since subgroups
defined by orders form an important class of examples, we do give
an explicit bound for the principal congruence subgroups of
$\Gamma$, depending on the volume of the cocompact quotient, as
well as some numerical characteristics of the order defining
$\Gamma$. Therefore, assume~$\Gamma = \Lat$ where~$\CQ$ is an
order in a division algebra~$\quat$ over~$K$, as before.

Now assume~$F = \R$, and let~$X_I = \dom{\Gamma(I)}{\HC^2}$
where~$I \normali \algint{K}$ is a given ideal. Let~$\mu$ denote
the hyperbolic measure on~$\HC^2$.

\begin{thm}
\label{51}
Let~$K$ be a totally real number field of dimension~$d$
over~$\,\,\Q$, and let~$D$ be a quaternion division algebra as in
\eqref{Ddef}. Let~$\Gamma = \Lat$ where~$\CQ$ is an order in~$D$.
%
% Let~$\kappa \in \algint {K}$ be such that~$\CQ \sub
% \frac{1}{\kappa}\O$ (where~$\O$ is defined in \eqref{standord}).
Let~$\nu = \mu(X_1)$ where~$X_1 = \dom {\Gamma} {\HC^2}$. Then the
surfaces~$X_I = \dom{\Gamma(I)}{\HC^2}$, where~$I \normali
\algint{K}$, satisfy
\[
\syspi_1(X_I) \geq \frac{4}{3}\left[\log(\genus(X_I)) -
        \left( \log(2^{3d-5} \pi^{-1}\nu \lam_{D,\CQ}) +
        o(1)\right) \right],
\]
where~$\lambda_{D,\CQ}$ is defined as in \eqref{lamdef}, ranging
over the primes of~$T_1\cup T_2$ which divide~$I$.
\end{thm}

Let~$R(D,\CQ)= 8^d \nu(\CQ) \lambda_{D,\CQ}$, where~$d=[K:\Q]$.
Theorem \ref{51} can be restated as follows.

\begin{cor}
We have
\[
\syspi_1(X_I) \geq \tfrac{4}{3}\log(\genus(X_I)) +
\tfrac{4}{3}\log(32 \pi) - \tfrac{4}{3} \log R(D,\CQ) + o(1),
\]
\end{cor}

Thus, finding a family of principal congruence subgroups with the best
systolic lower bound amounts to minimizing the expression~$R(D, \CQ)$
over all orders.

\begin{proof}[Proof of Theorem~\ref{51}]
Clearly~$X_I$ is a cover of degree~$\dimcol{\Gamma}{\Gamma(I)}$ of
$X_1$, hence we have from \Cref{3.2} that
\begin{eqnarray}\label{3.1}
4 \pi (\genus(X_I) - 1) & \leq & \mu(\dom{\Gamma(I)}{\HC^2})\nonumber \\
    & = & \dimcol{\Gamma}{\Gamma(I)} \cdot \nu \label{bm}\\
    & \leq & \nu \lam_{D,\CQ} \cdot \Norm(I)^{3}\nonumber
\end{eqnarray}
and so~$\Norm(I) \geq \(\frac{4\pi(g-1)}{\nu
\lam_{D,\CQ}}\)^{1/3}$ where~$\genus = \genus(X_I)$ is the genus.
Note that the first line in \eq{bm} is an equality if~$\Gamma(I)$
is torsion free. Let~$\pm 1 \neq x \in \Gamma(I)$. By
\eqref{boundell}
and \Tref{boundtr}, we have
\begin{eqnarray*}
\ell_x & > & 2 \log\(\abs{\Tr[{\M[2](\R)}]} (x)-1\) \\ & > & 2
    \log\(\frac{1}{2^{2d-2}} \Norm(I)^{2}-3\)\\
     & > & 2 \log\(\frac{1}{2^{2d-2}}
        \(\frac{4\pi(g-1)}{\nu \lam_{D,\CQ}}\)^{2/3}-3\) \\
    & > & 2 \log\(\frac{1}{2^{2d-2}}
        \(\frac{4\pi g}{\nu \lam_{D,\CQ}}\)^{2/3}\) + o(1)\\ %-\frac{4}{3}\epsilon \\
    & = & \frac{4}{3}\left[\log(g) -
        \(\log(2^{3d-5}\pi^{-1}\nu \lam_{D,\CQ}) + o(1)\)\right].
\end{eqnarray*}
%where~$\epsilon > 0$ is needed to balance the fourth inequality,
% and can be arbitrary small when~$g \ra \infty$.

The constant can be somewhat improved by taking the stronger version
of \Tref{boundtr} into account.
\end{proof}
%% For~$m$ odd,~$3d$ is replaced by~$3d/2$.
%%$$
%%\(1-\frac{1}{g}\)^{2/3}-\frac{3\cdot 2^{d-2}{\abs{\Norm(\kappa)}(\pi\nu)^{2/3}}}{g^{2/3}}
    % > e^{-2/3\epsilon}
%%$$

\Tref{12b} now follows immediately, since every arithmetic co-compact
lattice of~$\SL[2](\R)$ is by definition commensurable to one of the
lattices treated in \Tref{51}.

To prove \Tref{12d}, let~$F = \C$. Let~$K$ be a number field with one
complex embedding and~$d-2$ real ones, where~$d =
\dimcol{K}{\Q}$. As before, let~$D = (a,b)_{2,K}$ be a quaternion
division algebra over~$K$,~$\O$,~$\CQ$ and~$\kappa$ as in
\Tref{51}. Let~$\Gamma = \Lat$, and set~$X_1 =
\dom{\Gamma}{\HC^3}$.

\begin{thm}\label{52} % The complex case:~$F=\C$.
Let~$K$ be a number field of dimension~$d$ over~$\Q$ with a single
complex place, and let~$D$,~$\O$,~$\CQ$ and~$\kappa$ be as in
\Tref {51}. As before, let~$\Gamma = \Lat$. Let~$I_0 \normali
\algint{K}$ be an ideal such that~$\Gamma(I_0)$ is torsion free
(see \Cref{torf}), and set~$X_1 = \dom {\Gamma(I_0)}{\HC^3}$.

Then, the~$3$-manifolds~$X_I = \dom {\Gamma(I)} {\HC^3}$, where~$I
\sub I_0$, satisfy
\[
\syspi_1(X_I) \geq \frac{2}{3}\left[\log(\|X_I\|)) -
        \(\log(2^{3d-6} \|X_1\| \lam_{D,\CQ}) +
        o(1)\)\right],
\]
where~$\lambda_{D,\CQ}$ be defined as in \eqref{lamdef}, ranging
over the primes of~$T_1 \cup T_2$ which divide~$I$.
\end{thm}
\begin{proof}
We proceed as in the proof of \Tref{51}.  Since the simplicial
volume is multiplicative under covers where the covered space is
smooth, we have from \Cref{3.2} that
\begin{eqnarray*}
{}\|X_I\| & = & \dimcol{\Gamma(I_0)}{\Gamma(I)} \cdot \|X_1\|\\
    & \leq & \dimcol{\Gamma}{\Gamma(I)} \cdot \|X_1\|\\
    & \leq & \|X_1\| \lam_{D,\CQ} \cdot \Norm(I)^{3}
\end{eqnarray*}
and so~$\Norm(I) \geq \(\frac{\|X_I\|}{\|X_1\|
\lam_{D,\CQ}}\)^{1/3}$. Let~$\pm 1 \neq x \in \Gamma(I)$. By
\eqref{boundell} and \Tref{boundtr}, we have
\begin{eqnarray*}
\ell_x & > & 2 \log\(\abs{\Tr[{\M[2](\C)}]} (x)-1\) \\ & > & 2
    \log\(\frac{1}{2^{d-2}}
    \Norm(I)-3\)\\
     & > & 2 \log\(\frac{1}{2^{d-2}}
        \(\frac{\|X_I\|}{\|X_1\| \lam_{D,\CQ}}\)^{1/3}-3\) \\
    & > & 2 \log\(\frac{1}{2^{d-2}}
        \(\frac{\|X_I\|}{\|X_1\| \lam_{D,\CQ}}\)^{1/3}\) + o(1)\\ %-\frac{4}{3}\epsilon \\
    & = & \frac{2}{3}\left[\log(\|X_I\|) -
        \(\log(2^{3d-6} \|X_1\| \lam_{D,\CQ}) + o(1)\)\right],
\end{eqnarray*}
and again the constant may be improved by taking the stronger
version of \Tref{boundtr} into account.
\end{proof}

\Tref{12d} now follows, by the same argument as for \Tref{12b}.
Corollary~\ref{16} follows from the fact that for a closed hyperbolic
3-manifold~$M$, we have~$\| M \| = \frac{\vol(M)}{v_3}$ where~$v_3$ is
the volume of a regular ideal simplex in~$\HC^3$.

\section{The systole of Hurwitz surfaces}
\label{237}

In this section, we specialize the results of Section \ref{B} to the
lattice~$\Deltah$, defined as the even part of the group of
reflections in the sides of the~$(2,3,7)$ hyperbolic triangle.  We
follow the concrete realization of the~$\Z_2$-central
extension~$\tDeltah$ of~$\Deltah$ as the group of norm one elements in
an order of a quaternion algebra, given by N.~Elkies in
\cite[p.~39]{El0} and in \cite[Subsection~4.4]{El}.  The~$(2,3,7)$
case is is also considered in detail in~\cite[pp.~159-160]{MR}.

Let~$K$ denote the real subfield of~$\Q[\rho]$, where~$\rho$ is a
primitive \th{7} root of unity.  Thus~$K = \Q[\cc]$, where~$\cc =
\rho+\rho^{-1}$ has minimal polynomial~$\cc^3 + \cc^2 - 2\cc -1 =
0$. There are three embeddings of~$K$ into~$\R$, defined by
sending~$\cc$ to any of the numbers
\[
2\cos\(\tfrac{2\pi}{7}\), 2\cos\(\tfrac{4\pi}{7}\),
2\cos\(\tfrac{6\pi}{7}\).
\]
We view the first embedding as the `natural' one, and denote the
others by~$\s_1,\s_2 \co K \ra \R$.

Now let~$\quat$ denote the quaternion~$K$-algebra %
$$(\cc,\cc)_K = K[i,j \subjectto i^2 = j^2 = \cc,\, ji = -ij].$$ The
ring of integers of~$K$ is~$\algint{K} = \Z[\cc]$, and so this
presentation satisfies the condition of Lemma~\ref{2.1}.

\begin{prop}
\label{71}
The only two ramification places of~$D$ are the real embeddings~$\s_1$
and~$\s_2$.
\end{prop}

\begin{proof}
This fact is mentioned without proof in \cite[p.~96]{El}, and we
provide the easy argument for the sake of completeness.

The behavior of~$D$ over a real place is determined by the sign
of~$\eta$.  Since~$2 \cos(2\pi/7)$ is positive, the algebra~$D \tensor
\R$ is a matrix algebra.  On the other hand,~$\s_1(\eta) = \eta^2 - 2
= 2\cos(4\pi/7)$ and~$\s_2(\eta) = 2\cos(8\pi/7)$ are negative, and
so~$D \tensor_{\s_1} \R$ and~$D \tensor_{\s_2} \R$ are division
algebras.

Now let~$\fp$ be a prime of~$K$. If~$\fp$ is odd, then the norm
form
\[
\lam_1^2 - \eta \lam_2^2 - \eta \lam_3^2 + \eta^2 \lam_4^2
\]
is a four dimensional form over~$\algint{K}/\fp$, so it represents
zero non-trivially over~$\algint{K}/\fp$ \cite[Prop.~IV.4]{Serre}.  By
Hensel's lemma, such a representation can be lifted to a
representation over the completion~$K_{\fp}$, and so~$D \tensor
K_{\fp}$ has elements of zero norm, which are clearly zero divisors;
thus~$D \tensor K_{\fp} \isom \M[2](K_{\fp})$.

The same argument works for~$\fp = \ideal{2}$ (the only even
prime, noting that~$\algint{K}/\ideal{2} \isom \F_8$ is a field),
with one refinement necessary for Hensel's lemma to apply, namely,
the form represents zero modulo~$\fp^3 = 8$:
\[
1 - \eta(1+3\eta+\eta^2)^2 - \eta (\eta)^2 + \eta^2(0)^2 =
-8(1+3\eta+\eta^2) \equiv 0 \pmod{8}.
\]
As an alternative to the computation modulo~$8$, one can deduce the
behavior at~$\fp = 2$ from the other places, via the quadratic
reciprocity, which forces an even number of ramification places.
\end{proof}

% In fact by the Albert-Brauer-Hasse-Noether theorem,~$A$ is the
%  only quaternion algebra with this property (and the same behavior
%  with respect to the real embeddings).

Let~$\O \sub \quat$ be the order defined by \eqref{standord},
namely~$\O = \Z[\cc][i,j]$. Fix the element~$\tau = 1+\cc+\cc^2$, and
let
\[
j' = \frac{1}{2}(1+\cc i+\tau j).
\]
Notice that~$j'$ is an algebraic integer of~$\quat$, since the
reduced trace is~$1$ while the reduced norm is
\[
\frac{1}{4}(1-\cc\cdot \cc^2 - \cc\cdot \tau^2 + \cc^2 \cdot 0) =
-1-3\cc,
\]
so that both are in~$\algint{K}$. In particular
\begin{equation}\label{jt}
{j'}^2 = j' + (1+3\eta).
\end{equation}

\begin{defn}
We define the Hurwitz quaternion order~$\Hur$ by setting
\begin{equation}
\label{72}
\Hur = \Z[\cc][i,j,j'].
\end{equation}
\end{defn}

There is a discrepancy between the descriptions of a maximal order in
\cite[p.~39]{El0} and in \cite[Subsection~4.4]{El}.  The maximal order
according to \cite[p.~39]{El0} is~$\Z[\cc][i,j,j']$.  Meanwhile, in
\cite[Subsection~4.4]{El}, the maximal order is claimed to be the
order~$\Z[\cc][i,j']$, described as~$\Z[\cc]$-linear combinations of
the elements~$1$,~$i$,~$j'$, and~$ij'$ on the last line of page 94.
The correct answer is the former, \ie \eqref{72}.  More details may be
found in \cite{KSV2}.

Obviously~$\Hur \sub \frac{1}{2}\O$. Moreover since~$\Tr[\quat](j')
= 1$, this is the best possible choice for~$\kappa$.

\begin{lemma}
\label{72c}
In the Hurwitz case, one has~$\lambda_{D,\Hur}=1$.
\end{lemma}

\begin{proof}
Indeed, in~\eq{lamdef}, we have~$T_1 = \emptyset$ since~$D$ has no
finite ramification points by Proposition~\ref{71}, while~$T_2 =
\emptyset$ since the Hurwitz order~$\Hur$ is maximal.
\end{proof}

\begin{lemma}
\label{72b}
The subgroup~$\Hurone(I)$ is torsion-free for every proper ideal~$I
\normali \algint{K}$.
\end{lemma}

\begin{proof}
Let~$I \normali \algint{K}$ be a proper non-zero ideal, and assume
$\Hurone(I)$ is not torsion free. Taking into account the fact
that~$\algint{K}$ is a principal ideal domain, we have by
\Cref{torf} that~$I^2$ divides an ideal of the form
\[
\ideal{x^{-1}(x-1)^2} = \ideal{x+x^{-1}-2},
\]
where~$x$ a root of unity for which~$K[x]/K$ is a quadratic
extension. For our field~$K$, we must have~$x^{14} = 1$, namely~$x
= \pm \rho^j$,~$j = 1,\dots,6$. Now, the element
\[
2 + (\rho+\rho^{-1}) = 2 +\eta,
\]
as well as its Galois conjugates, is invertible in~$\Z[\eta]$ (as
$\eta(\eta-1)(\eta+2) = 1$), ruling out the cases~$x = - \rho^j$.
In all the remaining cases~$I = \ideal{2-\eta}$, since this ideal
is stable under the Galois action (as~$\s_1 \co 2-\eta\mapsto
(2-\eta)(2+\eta)$). However~$\ideal{2-\eta}$ is prime (of
norm~$7$), so cannot be divisible by~$I^2$.
\end{proof}

\begin{proof}[Proof of \Tref{special}]
For low genus, \eg~$\genus\leq 100$, one can verify directly that the
Hurwitz surfaces satisfy the~$4/3$-bound.  The Hurwitz surfaces with
automorphism group of order up to a million were classified by
M.~Conder~\cite{Co} using group theoretic arguments. The few surfaces
of genus below~$100$ can be dealt with on a case by case basis.  We
list them in Table~\ref{Hur}, together with their systoles and the
bound~$ \frac{4}{3} \log \genus(X)$. Here a surface is called {\em
simple\/} if its automorphism group is.  A surface is called {\em
real\/} if it admits an antiholomorphic involution.  Thus, if a
surface is non-real, then there are two distinct Riemann surfaces
which are isometric as Riemannian manifolds.  The systoles given in
the table were calculated by R.~Vogeler \cite[Appendix C]{Vo03}, \cf
\cite{Vo}.

%%%%
\begin{table}
\begin{tabular}{|r||l|l|l|l|l|l|}
\hline
 g  & Automorphism Group & Type & Reality  &Systole &Bound \\
\hline \hline
 3    &~$\PSL[](2,7)$ & simple & real  & 3.936 & 1.465 \\
 7    &~$\PSL[](2,8)$ & simple & real  & 5.796 & 2.595 \\
14    &~$\PSL[](2,13)$& simple & real  & 5.903  & 3.519 \\
14    &~$\PSL[](2,13)$& simple & real  & 6.887 & 3.519 \\
14    &~$\PSL[](2,13)$ & simple & real & 6.393 & 3.519 \\
17 &$(C_2)^3.\PSL[](2,7)$ & non-simple & non-real & 7.609 & 3.778 \\
\hline
%%% The next largest genus for which a Hurwitz surface exists
%%%  in~$118$.  It is a simple, real Hurwitz surface with automorphism
%%%  group~$PSL(2,27)$, systole 10.451 and bound 6.362.
\end{tabular}
\vskip 0.5cm
 \caption{Hurwitz surfaces of genus~$\leq 65$}
\label{Hur}
\end{table}
%%%%%
A comparison of the values in the last two columns, together with
Lemma~\ref{72b}, verifies the validity of the~$4/3$-bound for these
surfaces.  The general case is treated below.
\end{proof}

Theorem \ref{boundtr} specializes to the following result.

\begin{thm}
Let~$I \normali \Z[\eta]$.
For every~$x \neq \pm 1$ in~$\Deltah(I)$ we have
$$
  \abs{\Tr[D](x)} >
  \frac{1}{16}\Norm(I)^2 - 2.$$
\end{thm}
% We now continue as in Section \ref{finish}.

To complete the proof of Theorem~\ref{special}, we let~$X_I =
\dom{\Deltah(I)}{\HC^2}$. Combining this bound with \eqref{boundell}
and \eqref{3.1}, the fact that~$\nu = \pi/21$,
%~$\pi/21$ follows from Gauss-Bonnet theorem.
and applying Lemma~\ref{72c}, we obtain
%~$$\sys\pi_1(\dom{\Deltah(m)}{\HC^2}) > \frac{4}{3}
     % \( \log\(g\) - \(\log(\frac{2^{4+9\beta(m)/2}}{21} ) + o(1)\)\).$$
\begin{eqnarray*}
\syspi_1(X_I) & > &  2\log\( \(\frac{21 (g-1)}{16}\)^{2/3}  - 3\),
\end{eqnarray*}
which implies
$$\syspi_1(X_I) \geq \frac{4}{3} \log \genus(X_I)$$
if~$\genus(X_I) \geq 65$.

\section{Acknowledgments}
We are grateful to Misha Belolipetsky, Amiram Braun, Dima Kazhdan,
Alan Reid, David Saltman, and Roger Vogeler for helpful discussions.

\bibliographystyle{amsalpha}

\end{document}